\documentclass[a4paper,12pt]{amsart}

\usepackage{xcolor}
\usepackage{amsmath}
\usepackage{amsfonts}
\usepackage{amssymb}
\usepackage{amsthm}
\usepackage{amsopn}
\usepackage{amscd}
\usepackage[utf8]{inputenc}
\usepackage{courier}
\usepackage{url}
\usepackage{color}
\usepackage{paralist}

\usepackage[vcentering,dvips]{geometry}
\title[VersalDeformations]{VersalDeformations --- a package for computing versal deformations and local Hilbert schemes}

\author[N.O.~Ilten]{Nathan Owen Ilten}
\address{Department of Mathematics,
        University of California,
        Berkeley, 94720 California, USA}
\email{nilten@math.berkeley.edu}
\newcommand{\ut}{\hspace{.1cm}\underline{t}\hspace{.1cm}}
\newcommand{\PP}{\mathbb{P}}

\DeclareMathOperator{\spec}{Spec}
\DeclareMathOperator{\Hom}{Hom}
\DeclareMathOperator{\im}{Im}
\DeclareMathOperator{\jac}{Jac}
\DeclareMathOperator{\tr}{tr}

\theoremstyle{theorem}

\newtheorem*{theorem*}{Theorem}

\theoremstyle{definition}

\begin{document}
\maketitle
\begin{abstract}
We provide an overview of the Macaulay2 package \emph{VersalDeformations}, which algorithmically computes versal deformations of isolated singularities, as well as local (multi)graded Hilbert schemes.
\end{abstract}

\section{Introduction}
Deformation theory provides mathematicians with the tools to describe (local) parameter spaces for various algebraic geometric objects, for example, for isolated singularities, or for invertible sheaves on a projective variety, see \cite{Ha}. However, computing such spaces in practice can be quite difficult, and by hand often intractable. The Macaulay2 \cite{M2} package \emph{VersalDeformations}, available online \cite{def}, aims to facilitate such calculations for two concrete deformation problems: versal deformations of isolated singularities, and local (possibly multigraded) Hilbert schemes.

The package \emph{VersalDeformations} provides several functions that may be used to calculate tangent and obstruction spaces for the deformation problems mentioned above. The function \verb+normalMatrix+ may be used to calculate a basis for any degree of the normal module of some (multi)hom\-oge\-neous ideal in a polynomial ring. The scripted functor \verb+CT+ may be used to calculate bases of the first and second cotangent cohomology modules $T_A^1$ and $T_A^2$ of some algebra $A$ over a field $k$, assuming that these modules are finite dimensional vector spaces. In the (multi)homogeneous case, \verb+CT+ may also be used to calculate bases of homogeneous pieces of these modules.

The main contribution of the package is the method \verb+versalDeformation+, which uses the Massey product algorithm to iteratively lift solutions of a deformation equation to higher and higher order; we describe this in more detail in the following section. This can be used to find power series descriptions of versal deformations and local Hilbert schemes. Since such a description may not be polynomial, the package provides an interface allowing the user to control at what point the lifting should terminate. The package also implements a more time-consuming lifting algorithm (via the option \verb+SmartLift+) that seeks to minimize the number of higher order terms appearing in the equations for the parameter space.

There are a number of other software packages that provide related functionality. J.~Stevens has written scripts for the original Macaulay to calculate $T^1$ and $T^2$ for isolated singularities \cite{St}. There is a library for Singular \cite{Singular} by B.~Martin that calculates the versal deformation of an isolated singularity as well as of modules \cite{BM}. B.~Hovinen has written a package for Macaulay2 that computes versal deformations of maximal Cohen-Macaulay modules on hypersurfaces \cite{Ho}. Finally, J.~Boehm is developing a package for computations involving deformations of Stanley-Reisner rings \cite{SR}.

\section{Solving the Deformation Equation}
In the following, we briefly describe the Massey product algorithm as we have implemented it. For more details and mathematical background, see \cite{laudal:79a}, \cite{laudal:86a}, \cite{St}, or \cite{St2}. For simplicity, we restrict to the case of the versal deformation of an isolated singularity, although our approach for Hilbert schemes is similar.

First we fix some notation. Let $S$ be a polynomial ring over some field $k$, and let $I$ be an ideal of $S$ defining a scheme $X=\spec S/I$ with isolated singularity at the origin. Consider a free resolution of $S/I$:
$$
\begin{CD}
	\cdots @>>> S^l @>R^0>> S^m @>F^0>> S @>>> S/I @>>> 0.
\end{CD}
$$
Let $\phi_i\in\Hom(S^m/\im R^0,S)$ $i=1,\ldots,n$ represent a basis of 
$$
T_{S/I}^1\cong\Hom(S^m/\im R^0,S)/\jac F^0.
$$
We introduce deformation parameters $t_1,\ldots,t_n$ and consider the map
$F^1\colon S[\ut]^m\to S[\ut]$ defined as 
$$
F^1=F^0+\sum_{i=1}^n t_i\phi_i.
$$
Let $\mathfrak{m}$ be the ideal generated by $t_1,\ldots,t_n$.
It follows that there is a map $R^1\colon S[\ut]^l\to S[\ut]^m$ with 
$R^1\equiv R^0\mod\mathfrak{m}$ satisfying the first order deformation equation
$$
F^1R^1\equiv0\mod\mathfrak{m}^2.
$$

Our goal is to lift the above equation to higher order, that is, for each $i>0$, to find $F^i\colon S[\ut]^m\to S[\ut]$ with $F^i\equiv F^{i-1}\mod \mathfrak{m}^{i}$ and $R^i\colon S[\ut]^l\to S[\ut]^m$ with $R^i\equiv R^{i-1}\mod \mathfrak{m}^i$ satisfying $F^iR^i\equiv 0\mod\mathfrak{m}^{i+1}$. In general, there are obstructions to doing this, governed by the $d$-dimensional $k$ vector space $T_{S/I}^2$. Thus, we instead aim to solve
\begin{equation}\label{defeqn}
	(F^iR^i)^{\tr}+C^{i-2}G^{i-2}\equiv 0\mod\mathfrak{m}^{i+1},
\end{equation}
where 	$(F^iR^i)^{\tr}$ denotes the transpose of $(F^iR^i)$. 
Here, the matrices $G^{i-2}\colon k[\ut]\to k[\ut]^d$ and $C^{i-2}\colon S[\ut]^d\to S[\ut]^l$ are congruent modulo $\mathfrak{m}^i$ to $G^{i-3}$ and $C^{i-3}$, respectively. Furthermore, we require that $G^i$ and $C^i$ vanish for $i<0$, and $C^0$ is of the form $V\cdot D$, where $V\in\Hom(S^d,S^l)$ gives representatives of a basis for $T_{S/I}^2$ and $D\in\Hom(S^d,S^d)$ is a diagonal matrix.
The matrices $G^i$ now give equations for the miniversal base space of $X$.

Our implementation solves \eqref{defeqn} step by step.  Given a solution $(F^i,R^i,G^{i-2},C^{i-2})$ modulo $\mathfrak{m}^{i+1}$, the package uses Macaulay2's built in matrix quotients to first solve for $F^{i+1}$ and $G^{i-1}$ (by working over the ring $S[\ut]/I+\im (G^{i-2})^{\tr} +\mathfrak{m}^{i+2}$) and then solve for $R^{i+1}$ and $C^{i-1}$. For the actual computation, we avoid working over quotient rings involving high powers of $\mathfrak{m}$ by representing the $(F^i,R^i,G^{i-2},C^{i-2})$ as lists of matrices that keep track of the orders of the parameters $t_j$ involved.
\section{Examples}
We provide two examples: a versal deformation and a multigraded Hilbert scheme. We begin with the classical example of the miniversal deformation of the cone over the rational normal curve of degree $4$, see
\cite{minusfour}.
\begin{verbatim}
i1 : loadPackage "VersalDeformations";

i2 : S=QQ[x_0..x_4];

i3 : I=minors(2,matrix {{x_0,x_1,x_2,x_3},{x_1,x_2,x_3,x_4}});

o3 : Ideal of S

i4 : F0=gens I;

             1       6
o4 : Matrix S  <--- S

i5 : transpose F0

o5 = {-2} | -x_1^2+x_0x_2  |
     {-2} | -x_1x_2+x_0x_3 |
     {-2} | -x_2^2+x_1x_3  |
     {-2} | -x_1x_3+x_0x_4 |
     {-2} | -x_2x_3+x_1x_4 |
     {-2} | -x_3^2+x_2x_4  |

             6       1
o5 : Matrix S  <--- S
\end{verbatim}
We see that the tangent space $T_{S/I}^1$ of the miniversal deformation is four-dimensional, and the obstruction space $T_{S/I}^2$ is  three-dimensional: 
\begin{verbatim}
i6 : CT^1(F0)

o6 = {-2} | x_1  x_0  0   0    |
     {-2} | 0    0    0   x_0  |
     {-2} | -x_3 -x_2 0   x_1  |
     {-2} | 0    0    x_2 0    |
     {-2} | -x_4 -x_3 x_3 0    |
     {-2} | 0    0    x_4 -x_3 |

             6       4
o6 : Matrix S  <--- S

i7 : CT^2(F0)

o7 = {-3} | 0    0   0   |
     {-3} | 0    0   0   |
     {-3} | x_3  x_4 0   |
     {-3} | x_2  x_3 0   |
     {-3} | x_1  x_2 0   |
     {-3} | -x_4 0   x_3 |
     {-3} | 0    x_4 x_2 |
     {-3} | 0    x_3 x_1 |

             8       3
o7 : Matrix S  <--- S
\end{verbatim}
In this example, our algorithm gives a polynomial solution to the deformation equation:
\begin{verbatim}
i8 : (F,R,G,C)=versalDeformation(F0,Verbose=>2);
Calculating first order deformations and obstruction space
Calculating first order relations
Starting lifting
Order 2
Order 3
Solution is polynomial
\end{verbatim}
The miniversal base space is the union of $\mathbb{A}^3$ and $\mathbb{A}^1$, intersecting in a point, as can be seen from the following equations:
{\small\begin{verbatim}
i9 : T=ring first G;

i10 : sum G

o10 = | t_2t_3-t_3^2 |
      | t_1t_3       |
      | t_3t_4       |

              3       1
o10 : Matrix T  <--- T
\end{verbatim}}
A versal family is given by
{\small\begin{verbatim}
i11 : transpose sum F

o11 = {0, -2} | t_1x_1+t_2x_0-x_1^2+x_0x_2                |
      {0, -2} | t_4x_0-x_1x_2+x_0x_3                      |
      {0, -2} | -t_1t_4-t_1x_3-t_2x_2+t_4x_1-x_2^2+x_1x_3 |
      {0, -2} | t_2t_3-t_3^2+t_3x_2-x_1x_3+x_0x_4         |
      {0, -2} | t_3t_4-t_1x_4-t_2x_3+t_3x_3-x_2x_3+x_1x_4 |
      {0, -2} | t_3x_4-t_4x_3-x_3^2+x_2x_4                |

              6       1
o11 : Matrix T  <--- T
\end{verbatim}}

We now consider our second example: the local description of the Hilbert scheme of the diagonal in $\PP^2\times\PP^2\times\PP^2$ at the point corresponding to the unique Borel fixed ideal, see  \cite{multigraded} for more details.
{\small\begin{verbatim}
i12 : S=QQ[x1,x2,x3,y1,y2,y3,z1,z2,z3,Degrees=>{{1,0,0},{1,0,0},
      	  {1,0,0},{0,1,0},{0,1,0},{0,1,0},{0,0,1},{0,0,1},{0,0,1}}];

i13 : I=ideal {y1*z2, x1*z2, y2*z1, y1*z1, x2*z1, x1*z1, x1*y2, x2*y1,
      	   x1*y1, x2*y2*z2};

o13 : Ideal of S

i14 : (F,R,G,C)=versalDeformation(gens I,normalMatrix({0,0,0},gens I),
           CT^2({0,0,0},gens I),Verbose=>2);
Calculating first order relations
Starting lifting
Order 2
Order 3
Order 4
Order 5
Order 6
Solution is polynomial
\end{verbatim}}

Note that since we were interested in the multigraded Hilbert scheme, the tangent space is just the degree $(0,0,0)$ component of the normal module of $I$, and an obstruction space is given by the degree $(0,0,0)$ component of $T_{S/I}^2$. In any case, this multigraded Hilbert scheme is locally cut out by 8 cubics:
{\small\begin{verbatim}
i15 : T=ring first G;

i16 : sum G

o16 = | t_2t_3t_4-t_2t_4t_7-t_1t_3t_8+t_1t_7t_8+t_1t_3t_13-...
      | t_1t_3t_4-t_2t_3t_4-t_1t_7t_8+t_2t_7t_8-t_1t_3t_13+...
      | t_1t_3t_16-t_2t_7t_16-t_1t_14t_16+t_2t_14t_16-...
      | t_1t_3t_18-t_2t_7t_18-t_1t_14t_18+t_2t_14t_18-...
      | t_2t_4t_17-t_1t_8t_17+t_1t_13t_17-t_2t_13t_17-...
      | t_2t_4t_18-t_1t_8t_18+t_1t_13t_18-t_2t_13t_18-...
      | t_3t_4t_17-t_7t_8t_17-t_3t_13t_17+t_7t_13t_17-...
      | t_3t_4t_16-t_7t_8t_16-t_3t_13t_16+t_7t_13t_16-...
\end{verbatim}}
There are in fact $7$ irreducible components of the Hilbert scheme that pass through this point:
{\small\begin{verbatim}
i17 : # primaryDecomposition ideal sum G

o17 = 7
\end{verbatim}}

\bibliographystyle{alpha}
\bibliography{versal}

\end{document}